\input amstex
\documentstyle{amsppt}
\magnification=\magstep1

\define\sa{\operatorname{sa}}
\define\Tr{\operatorname{\,Tr}}
\define\diag{\operatorname{\,diag\,}}
\define\ua{\underline{a}}
\define\ub{\underline{b}}
\redefine\i{\operatorname{\,i}}
\redefine\sp{\operatorname{\,sp}}
\overfullrule=0pt
\TagsOnRight

\topmatter  
\title{Jensen's operator inequality}\endtitle 

\rightheadtext{Jensen's Operator Inequality}

\author Frank Hansen \& Gert K. Pedersen \endauthor

\date  {$2^{nd}$ April, 2002}\enddate

\address {Department of Economics, University of Copenhagen,
Studiestr\ae{}de 6, DK-1455 Copenhagen K, Denmark
\& Department of Mathematics,  University of Copenhagen, 
Universitetsparken 5, DK-2100 Copenhagen \O, Denmark} 
\endaddress

\email {frank.hansen\@ econ.ku.dk \& gkped\@ math.ku.dk}\endemail

\abstract{We establish what we consider to be the definitive
versions of Jensen's operator inequality and Jensen's 
trace inequality for real functions defined on an 
interval. This is accomplished by the introduction of 
genuine non-commutative convex combinations of operators, 
as opposed to the contractions considered in earlier 
versions of the theory, \cite{{\bf 9}} \& \cite{{\bf 3}}.
As a consequence, we no longer need to impose 
conditions on the interval of definition. We  
show how this relates to the pinching inequality of 
Davis \cite{{\bf 4}}, and how Jensen's  trace 
inequality generalizes to $ C^*-$algebras.}
\endabstract

\endtopmatter
\document       
\footnote""{\copyright 2002 by the authors. This paper may be 
reproduced, in its entirety, for non-commercial purposes.} 

\head{1. Introduction}\endhead

If $f$ is a continuous, real function on some interval 
$I$ in $\Bbb R$, we can use spectral theory to define 
an operator function
$$
f\colon \Bbb B(\frak H)^I_{\sa} @>>> \Bbb B (\frak H)_{\sa} \quad
\text{where}\quad f(x)=\int f(\lambda)\,dE_x(\lambda).
\tag{1}
$$
Here $\Bbb B(\frak H)^I_{\sa}$ denotes the convex set 
of self-adjoint operators on the Hilbert space $\frak H$ 
with spectra in $I$, and $E_x$ denotes the spectral 
measure of $x$. Admittedly it is somewhat dangerous to 
use the same symbol for the two rather different functions, 
but the usage is sanctified by time. Whenever necessary 
we shall try to distinguish between the two by referring 
either to the {\it function} $f$ or to the 
{\it operator function} $f$. As pointed out by C. Davis 
in \cite{{\bf 4}} a general operator function 
$F\colon\Bbb B(\frak H)^I_{\sa}@>>>\Bbb B(\frak H)$ 
arises from a spectral function, i\.e\. $F(x)=f(x)$, 
if and only if for every unitary $u$ 
$$
F(u^*xu)=u^*F(x)u \quad \text{and}\quad 
F\left(\matrix y&0\\0&z\endmatrix\right) 
= \left(\matrix F(y) & 0\\0 &
F(z)\endmatrix\right)
\tag{2}
$$
for every operator $x=y+z$ that decomposes in block 
form by multiplication by a projection $p$ in 
its commutant. (N\.b\. we do not demand that 
$p$ and $\bold 1-p$ are equivalent.) There is a 
slight ambiguity in this statement -- easily 
compensated for by its versatility -- since by 
F(y) we really mean $F$ evaluated at $y$, but 
now regarded as an operator function on 
$\Bbb B(p\frak H)^{I}_{\sa}$. Put differently, 
we demand that $pF(y+z)=F(y+z)p$ and that it is 
independent of $z$. Thus, $pF(y+z)= 
pF(y+s(\bold 1 -p))$ for some, hence any scalar 
$s$ in $I$. (Davis tacitly assumes that $0\in I$ 
and takes $s=0$.)

A continuous function $ f\colon I @>>> \Bbb R$ 
is said to be {\it operator convex} if
$$
f(\lambda x+(1-\lambda)y)\le\lambda f(x)+(1-\lambda)f(y)
\tag{3}
$$
for each $\lambda$ in $[0,1] $ and every pair of 
self-adjoint operators $x, y$ on an infinite dimensional 
Hilbert space $\frak H$ with spectra in $I$.
The function is said to be {\it matrix convex of order 
$n$} if the same conditions are satisfied for operators 
on a Hilbert space of finite dimension $n$. 
It is well known, cf\. \cite {{\bf 2}, Lemma 2.2} that 
a function is operator convex if and only if it is 
matrix convex of arbitrary orders. 

Just because the function $f$ is convex there is no 
guarantee that the operator function $f$ is convex. 
In fact, as shown by Bendat and Sherman in \cite{{\bf 2}}, 
$f$ is operator convex on the interval $]-1, 1[$ if 
and only if it has a (unique) representation 
$$
f(t)= \beta_0 + \beta_1t + \tfrac 12 \beta_2\int_{-1}^1 t^2(1-\alpha
t)^{-1}\,d\mu(\alpha),
\tag{4}
$$
for $\beta_2 \ge 0$ and some probability measure 
$\mu$ on $[-1, 1]$. In particular, $f$ must be analytic 
with $f(0)=\beta_0, \; f'(0)=\beta_1$ and $f''(0)=\beta_2$. A 
concise account of this result and its relations to
L\"o{}wner's theory of operator monotonicity can  be 
found in \cite{{\bf 9}}. 

An unexpected phenomenon turns up in relation with 
convexity in $\Bbb B(\frak H)_{\sa}$.  If 
$(a_1, \dots, a_n)$ is an $n-$tuple of operators with
$\sum_{k=1}^n a^*_ka_k=\bold 1$, we may think of the 
element  $\sum_{k=1}^n a^*_kx_ka_k$ as a non-commutative 
convex combination of the  $n-$tuple $(x_1, \dots, x_n)$ 
in $\Bbb B(\frak H)_{\sa}$. The remarkable  fact is that
when $f$ is an operator convex function, then the 
operator function $f$ respects this new structure in the 
sense that we have the  {\it Jensen operator inequality}:
$$
f\left(\sum_{k=1}^n a^*_kx_ka_k\right) 
\le \sum_{k=1}^n a^*_kf(x_k)a_k.
\tag{5}
$$
This result was found in embryonic form by the first 
author in \cite {{\bf 6}}, and used by the two of us 
to give a review of L\"o{}wner's and Bendat-Sherman's 
theory of operator monotone and operator convex 
functions in \cite {{\bf 9}}.  With hindsight we must 
admit that we unfortunately proved and used the 
contractive form $f(a^*xa)\le a^*f(x)a$ for
$a^*a\le\bold 1$,  this being the seemingly most 
attractive version at the time. However, this 
necessitated the further conditions that $0\in I$ and
$f(0)\le 0$, conditions that have haunted the theory 
since then. The Jensen inequality for a normal trace 
on a von Neumann algebra, now for an arbitrary convex 
function $f$, was found by Brown and Kosaki in 
\cite{{\bf 3}}, still in the contractive version.

It is the aim of the present paper to rectify our 
omissions and prove the full Jensen inequality, 
both with and without a trace.  This is accomplished 
by a refinement of previous techniques and by 
applying some new ideas that also make the 
presentation more streamlined and easier to follow.

\vskip 2truecm

\head{2. Main Results}\endhead


\proclaim{Theorem 2.1 {\rm (Jensen's Operator Inequality)}}

\noindent For a continuous function $f$ defined on an 
interval $I$ the following conditions are equivalent:

\smallskip

\noindent{\rm (i)}\; $f$ is operator convex. 

\smallskip

\noindent{\rm (ii)}\; For each natural number $n$ we 
have the inequality
$$
f\left(\sum_{i=1}^n a_i^* x_i a_i\right)
\le\sum_{i=1}^n a_i^* f(x_i)a_i
\tag{5}
$$
for every $n-$tuple $(x_1,\dots, x_n)$ of bounded, 
self-adjoint operators on an arbitrary Hilbert 
space $\frak H$ with spectra contained in $I$ and 
every $n-$tuple $(a_1,\dots, a_n)$ of operators on 
$\frak H$ with $\sum_{k=1}^n a_k^*a_k = \bold 1$.

\smallskip

\noindent{\rm (iii)}\; $f(v^*xv)\le v^*f(x)v$  for each 
isometry $v$ on an infinite-dimensional Hilbert space  
$\frak H$ and every self-adjoint operator  $x$ with 
spectrum in $I$.

\smallskip

\noindent{\rm (vi)}\; $pf(pxp+s(1-p))p \le pf(x)p$ for 
each projection $p$ on an infinite-dimensional 
Hilbert space $\frak H$, every self-adjoint operator 
$x$ with spectrum in $I$ and every $s$ in $I$.
\endproclaim


\subhead{Remark 2.2}\endsubhead If the Jensen operator 
inequality (5) is satisfied for some $n\ge 2$ and 
for operators on a Hilbert space $\frak H$ 
(of any dimension), then clearly $f$ is operator
convex on $\Bbb B(\frak H)$. The point of condition (iii)  
is that if $\frak H$ infinite-dimensional it suffices to take 
$n=1$.  On the other hand it is clear that if (5) is satisfied 
for some $n$, then (setting $a_i=0$ for $i>1$) it is also 
satisfied for $n=1$.


\proclaim{Corollary 2.3 {\rm (Contractive Version)}}

\noindent Let $f$ be a continuous function defined on 
an interval $I$ and suppose that $0\in I$. Then $f$ is 
operator convex and $f(0)\le 0$ if and only if for 
some, hence every natural number $n$, the inequality 
{\rm (5)} is valid for every $n-$tuple $(x_1,\dots, x_n)$ 
of bounded, self-adjoint operators on a Hilbert space 
$\frak H$ with spectra contained in $I$, and
every $n-$tuple $(a_1,\dots, a_n)$ of operators on 
$\frak H$ with $\sum_{k=1}^n a^*_ka_k \le \bold 1$.
\endproclaim

Setting $n=1$ we see that $f$ is operator convex on an 
interval $I$ containing $0$ with $f(0)\le 0$ if and only if 
$$
f(a^*xa)\le a^*f(x)a 
\tag{6}
$$
for every self-adjoint $x$ with spectrum in $I$ and every 
contraction $a$. This is the original Jensen operator 
inequality from \cite{\bf 9}.


\proclaim{Theorem 2.4 {\rm (Jensen's Trace Inequality)}} 

\noindent Let $f$ be a continuous function defined on an 
interval $I$ and let $m$ and $n$ be natural numbers. If 
$f$ is convex we then have the inequality 
$$
\Tr \left(f\left(\sum_{i=1}^n a_i^*x_ia_i\right)\right)
\le\Tr\left(\sum_{i=1}^n a_i^*f(x_i)a_i\right)
\tag{7}
$$
for every $n-$tuple $(x_1,\dots, x_n)$ of self-adjoint 
$m\times m$ matrices with spectra contained in $I$ and 
every $n-$tuple $(a_1,\dots, a_n)$ of $m\times m$ matrices 
with $\sum_{k=1}^n a_k^*a_k=\bold 1$.

Conversely, if the inequality {\rm (7)} is satisfied for some $n$ and 
$m$, where $n>1$, then $f$ is convex.   
\endproclaim 


\proclaim{Corollary 2.5 {\rm (Contractive Version)}} 

\noindent Let $f$ be a convex, continuous function 
defined on an interval $I$, and suppose that 
$0\in I$ and $f(0)\le 0$. Then for all natural numbers 
$m$ and $n$ we have the inequality {\rm(7)} for every 
$n-$tuple $(x_1,\dots, x_n)$ of self-adjoint 
$m\times m$ matrices with spectra contained in $I$ and 
every $n-$tuple $(a_1,\dots, a_n)$ of $m\times m$ 
matrices with $\sum_{k=1}^n a_k^*a_k\le\bold 1$.  
\endproclaim


\subhead{Remark 2.6}\endsubhead Let $n=1$ in (7).  
If $f$ is convex, $0\in I$ and $f(0)\le 0$ we have
$$
\Tr\left(f(a^*xa)\right)\le \Tr\left(a^*f(x)a \right)
\tag{8}
$$
for every self-adjoint $m\times m$ matrix $x$ with 
spectrum in $I$ and every $m\times m$ contractive 
matrix $a$. This is Jensen's trace inequality (for 
matrices) of Brown and Kosaki [{\bf 3}]. 

This inequality alone is not sufficient to ensure 
convexity of $f$, even if $m>1$ (unless $f(0)=0$ is 
specified in advance).  However, for $n>1$ the inequality 
gives convexity of $f$ as we see from Theorem 2.4. 
In each case we must assume that $0\in I$, otherwise the 
inequality does not make sense. This fact, together with the 
irrelevant information about $f(0)$,  makes the contractive 
versions of Jensen's inequality less desirable. When we 
eventually pass to the theory of several variables, 
cf\. \cite{{\bf 10}}, the contractive versions mean that 
$0$ belongs to the cube where $f$ is defined, so that 
part of the coordinate axes must belong to the domain 
of definition for $f$, and on these we must assume that 
$f\le 0$. This assumption is so severe a restraint that it 
becomes a real problem for the theory.


\proclaim{Theorem 2.7 {\rm (Jensen's Trace Inequality for 
$C^*-$Algebras)}}

\noindent Let $f$ be a convex, continuous function 
defined on an interval $I$ and let $\Cal A$ be a 
$C^*-$algebra with a finite trace $\tau$. Then the 
inequality 
$$
\tau\left(f\left(\sum_{i=1}^n
a_i^*x_ia_i\right)\right)\le
\tau\left(\sum_{i=1}^n a_i^*f(x_i)a_i\right)
\tag{9}
$$
is valid for every $n-$tuple $(x_1,\dots, x_n)$ of 
self-adjoint elements in $\Cal A$ with spectra    
contained in $I $ and every $n-$tuple 
$(a_1,\dots, a_n)$ in $ \Cal A$ with 
$\sum_{k=1}^n a^*_ka_k=\bold 1$.

If the trace $\tau$ is unbounded, but lower 
semi-continuous and densely defined, the 
inequality {\rm (9)} is still valid if $f\ge 0$, 
although now some of the numbers may be infinite. 
\endproclaim

\vskip 2truecm

\head{3. Unital and Unitary Tuples}\endhead 

\subhead{Notations}\endsubhead An $n-$tuple 
$\ua =(a_1, \dots , a_n)$ of operators in 
$\Bbb B(\frak H)$ is called a {\it contractive column}
(respectively a {\it unital column}) if 
$\sum_{k=1}^n a^*_ka_k \le \bold 1$ (respectively  
$\sum_{k=1}^n a^*_ka_k = \bold 1$). Contractive rows 
and unital rows are defined analogously by the 
conditions $\sum_{k=1}^n a_ka^*_k \le\bold 1$ and 
$\sum_{k=1}^n a_ka^*_k = \bold 1$. We say that 
$\ua = (a_1,\dots , a_n)$ is a {\it unitary column} 
if there is a unitary $n\times n$ operator matrix 
$U=(u_{ij})$, one of whose columns is 
$(a_1,\dots , a_n)$. Thus, $u_{ij}=a_i$ for some $j$ 
and all $i$. Unitary rows are defined analogously, 
cf\. \cite{{\bf 1}, Definition 1.1.} Note that an
$n-$tuple $(a_1,\dots, a_n)$ is a 
contractive/unital/unitary row if and only if the 
adjoint tuple $(a_1^*,\dots, a_n^*)$ is a
contractive/unital/unitary column.
Even for a finite-dimensional Hilbert space 
$\frak H$ it may happen that an $n-$tuple $\ua$ is 
a unitary (unital or contractive) column in 
$\Bbb B(\frak H)$, while $\ua$ is not a unitary
(unital or contractive) row in $\Bbb B(\frak H)$.
Evidently every unitary column is also a unital 
column (and similarly every unitary row is a unital 
row). On the other hand, if $(s_1, \dots , s_n)$
is an $n-$tuple of co-isometries such that 
$\sum_{k=1}^n s^*_ks_k =\bold 1$ (these are the 
canonical generators for the Cuntz algebra 
$\Cal O_n$), then we have a simple example of a 
unital column that is not unitary. If we insist 
that a unital column of elements in a unital 
$C^*-$algebra $\Cal A$ should be called a unitary 
column only if we can choose the unitary in 
$\Bbb M_n(\Cal A)$, then already for $\Cal A$
commutative, viz\.  $\Cal A=C(\Bbb S^5)$,  we have 
a unital $3-$column that is not a unitary column 
in $\Bbb M_3(\Cal A)$, cf\. \cite{{\bf 17}, 
Example 14}. 

Given a unital column $(a_1,\dots, a_n)$ we may 
regard it as an isometry 
$\ua\colon\frak H @>>> \frak H^n$, where 
$\frak H^n = \oplus_{i=1}^n \frak H$. Better still 
we may regard it as a partial isometry 
$V\colon \frak H^n @>>>\frak H^n$, where 
$V|\frak H^{n-1} =0$. Evidently the column is 
unitary precisely when $V$ extends to a 
unitary operator on $\frak H^n$, and this happens 
if and only if the index of $V$ is $0$, in the 
generalized sense that 
$\dim\ker V^* = (n-1)\dim \frak H$. Here 
$V^*(\xi_1,\dots,\xi_n)= a^*_1\xi_1+\cdots + a^*_n\xi_n$ 
in $\frak H$. It follows from \cite{{\bf 1}, 
Corollary 2.2} that this holds if just one 
of the operators $a_i$ has (generalized) index zero, 
since in this case $a_i=u|a_i|$ for some unitary $u$ 
on $\frak H$. We are then reduced to the situation 
where one of the operators, say $a_n$, is positive, 
so that with $\ub = (a_1,\dots, a_{n-1})$ we can extend $V$ 
to the unitary operator
$$
U=\pmatrix (\bold 1 -\ub(\ub)^*)^{1/2} & \ub\\
 -(\ub)^* & a_n \\ \endpmatrix.
\tag{10}
$$
It follows that every contractive $n-$column can be 
enlarged to a unitary $(n+1)$--column simply by 
setting $a_{n+1}=(\bold 1-\sum_{k=1}^n a^*_ka_k)^{1/2}$. 
In particular, every unital $n-$column can be enlarged 
to a unitary $(n+1)-$column with $a_{n+1}=0$. As usual 
we shall refer to this as a {\it unitary dilation} of 
the unital (or contractive) column. 


\subhead{Unitary Dilations}\endsubhead It may 
sometimes be desirable to know exactly the terms in 
a unitary dilation of some unital column 
$\ua= (a_1, \dots , a_n)$. If $n=1$, so that $\ua=a$ for 
some isometry $a$, the canonical dilation is given by a
$2\times 2-$matrix $U$ having $(a, 0)$ as the second 
column. For a general unital $n-$column we may regard 
it as an isometry $\ua\colon\frak H\to\frak H^n$, and
the unitary dilation $U_n$ on $\frak H \oplus \frak H^n$ 
then has the same form as $U$; in fact 
$$
U=\pmatrix 1-aa^* & a
\\ -a^* & 0\endpmatrix \quad \text{and}\quad U_n=
\pmatrix p&\ua\\-(\ua)^*&0\endpmatrix,
\tag{11}
$$
where $p=\bold 1 - \ua(\ua)^*$ is the $n\times n$ 
projection in $\frak H^n$ with $p_{ii}= \bold 1 - a_ia_i^*$ and 
$p_{ij}= -a_ia^*_j$ for $i\ne j$. Thus, the {\it canonical dilation} 
of $(a_1, \dots , a_n)$ has the form:
$$
U_n = 
\pmatrix \bold 1 - a_1a_1^* & \;\;\;- a_1a_2^* & \dots & \;\;\; - a_1a_m^* 
& a_1\\
\;\;\; - a_2a_1^*& \bold 1 - a_2a_2^* & \dots & \;\;\; - a_2a_n^* & a_2\\
\vdots & \vdots & {} & \vdots & \vdots\\
\;\;\;  - a_na_1^* & \;\;\; - a_na_2^* & \dots & \bold 1 -a_na_n^* & a_n\\
\;\;\;\;\;\;\; - a_1^* & \;\;\;\;\;\;\; - a_2^* & \dots &
\;\;\;\;\;\;\; -a_n^* &  0  \endpmatrix 
\tag{12}
$$
As seen from (10), the formula for the canonical 
dilation of a contractive column is only marginally 
more complicated, cf\. \cite{{\bf 1}, Lemma 1.1}. 


\proclaim{Lemma 3.1} Define the unitary matrix 
$E =\diag(\theta,\theta^2,\dots,\theta^{n-1},1)$ in 
$\Bbb M_n(\Bbb C) \subset \Bbb B(\frak H^n)$, 
where $\theta = \exp\, (2\pi \i/ n)$.  Then for each element 
$A=(a_{ij})$ in $\Bbb B(\frak H^n)$ we have
$$
\frac 1n \sum_{k=1}^n E^{-k}AE^k 
= \diag(a_{11}, a_{22}, \dots , a_{nn}).
\tag{13}
$$
\endproclaim

\demo{Proof} By computation
$$
\left(\frac 1n \sum_{k=1}^n E^{-k}AE^k\right)_{ij} = \frac 1n \sum_{k=1}^n
(\theta^{j-i})^ka_{ij},
\tag{14}
$$
and this sum is zero if $i\ne j$, otherwise it is $a_{ii}$.
\hfill$\square$
\enddemo


\proclaim{Corollary 3.2} Let $P$ denote the projection 
in $\Bbb M_n (\Bbb C)$ given by $P_{ij}= n^{-1}$ for all 
$i$ and $j$, so that $P$ is the projection of rank one 
on the subspace spanned by the vector $\xi_1+\cdots +
\xi_n$ in $\Bbb C^n$, where $\xi_1,\dots, \xi_n$ are 
the standard basis vectors. Then with $E$ as in 
Lemma 3.1 we obtain the pairwise orthogonal projections 
$P_k= E^{-k}PE^k$, for $1\le k\le n$, with 
$\sum_{k=1}^n P_k = \bold 1$. 
\endproclaim

\demo{Proof} (Cf\. \cite{{\bf 8}, Proof of Theorem 2.6}) 
Evidently each $P_k$ is a projection of 
rank one. Moreover, by Lemma 3.1,
$$
\sum_{k=1}^n P_k = \sum_{k=1}^n E^{-k}PE^k 
= n\diag(n^{-1},\dots, n^{-1})=\bold 1,
\tag{15}
$$
from which it also follows that the projections are 
pairwise orthogonal,
\hfill$\square$
\enddemo

\vskip 2truecm

\head{4. Proofs and Further Results}\endhead

\demo{Proof of Theorem 2.1} (i)$\implies$(ii) Assume that 
we are given a unitary $n-$column $(a_1,\dots, a_n)$, and 
choose a unitary $U_n=(u_{ij})$ in $\Bbb B (\frak H^n)$ 
such that $u_{kn}=a_k$. Let $E =\diag (\theta, \theta^2,\dots, 1)$ 
as in Lemma 3.1 and put $X=\diag(x_1, \dots , x_n)$, both 
regarded as elements in $\Bbb B(\frak H^n)$. Using Lemma 3.1 
and the operator convexity of $f$ we then get the 
desired inequality:
$$
\aligned
f&\left(\sum_{k=1}^n a^*_kx_ka_k\right) 
= f\left((U_n^*XU_n)_{nn}\right)\\
=f&\left(\left(\sum_{k=1}^n \frac 1n E^{-k}U_n^*XU_nE^k\right)_{nn}\right)
=\left(f\left(\sum_{k=1}^n \frac 1n E^{-k}U_n^*XU_nE^k\right)\right)_{nn}\\
\le\;&\left(\frac 1n \sum_{k=1}^n f(E^{-k}U_n^*XU_nE^k)\right)_{nn}
=\left(\frac 1n \sum_{k=1}^n E^{-k}U_n^*f(X)U_nE^k\right)_{nn}\\
=\;&\left(U_n^*f(X)U_n \right)_{nn}=\sum_{k=1}^n a_k^*f(x_k)a_k.
\endaligned
\tag{16}
$$
Note that for the second equality we use that 
$f(y_n)= \left(f(\diag(y_1,\dots, y_n))\right)_{nn}$ because
$f(\diag(y_1,\dots, y_n)) = \diag(f(y_1),\dots, f(y_n))$.
 
In the general case where the column is just unital, 
we enlarge it to the unitary $(n+1)-$column 
$(a_1, \dots, a_n, 0)$ and choose $x_{n+1}$ arbitrarily, 
but with spectrum in $I$. By the first part of the proof 
we therefore have
$$
\aligned
&f \left(\sum_{k=1}^n a^*_kx_ka_k \right) 
=f\left(\sum_{k=1}^{n+1} a^*_kx_ka_k\right)\\
&\;\le\sum_{k=1}^{n+1} a^*_kf(x_k)a_k 
= \sum_{k=1}^n a^*_kf(x_k)a_k.
\endaligned
\tag{17}
$$

\medskip

\noindent(ii)$\implies$(iii) is trivial.

\medskip

\noindent(iii)$\implies$(iv) Take any self-adjoint 
operator $x$ with spectrum in $I$ and let $p$ be an 
infinite-dimensional projection. Then we can find 
an isometry $v$ (i\.e\. $v^*v=\bold 1$) such that  
$p=vv^*$. By assumption $f(v^*xv) \le v^*f(x)v$, 
whence also
$$
vf(v^*xv)v^* \le vv^*f(x)vv^*=pf(x)p.
\tag{18}
$$
For any monomial $g(t)=t^m$ and any $s$ in $I$ we have 
$$
\aligned
pvg(v^*xv)v^*p=\;&pv(v^*xv)^mv^*p = p(vv^*xvv^*)^m p\\ 
= pg(pxp)p =\;& pg(pxp+s(\bold 1 -p))p.
\endaligned
\tag{19}
$$
Since $f$ is continuous, it can be approximated by 
polynomials on compact subsets of $I$, and therefore 
also $pvf(v^*xv)v^*p=pf(pxp+s(\bold 1-p))p$. Combined 
with (17) this gives the pinching inequality 
$$
pf(pxp+s(\bold 1 -p))p \le pf(x)p.
\tag{20}
$$
If  $p$ is a projection of finite rank we can define 
the infinite dimensional projection 
$q=p\otimes \bold 1_{\infty}$ on $\frak H^{\infty}$. 
Similarly we let $y=x\otimes \bold 1_{\infty}$ for any
given self-adjoint operator $x$ with spectrum in $I$. 
Since $f(a)\otimes \bold 1_{\infty}= 
f(a\otimes \bold 1_{\infty})$ for any operator $a$ we 
get by (20) that
$$
\aligned
pf(pxp+s(\bold 1 -p))p \otimes \bold 1_{\infty} 
=\;&  qf(qyq+s(\bold 1-q))q\\
\le qf(y)q=\;&pf(x)p \otimes \bold 1_{\infty},
\endaligned
\tag{21}
$$
which shows that (20) is valid also for projections of 
finite rank.

\medskip

\noindent (iv)$\implies$(1)  Given self-adjoint 
operators $x$ and $y$ with spectra in $I$ and 
$\lambda$ in $[0,\, 1]$, define the three elements 
$$
X=\left(\matrix x &0\\0&y\endmatrix\right), 
\quad U=\left(\matrix\lambda^{1/2} & (1-\lambda)^{1/2}\\
 -(1-\lambda)^{1/2} &\lambda^{1/2}\endmatrix\right), 
\quad P=\left(\matrix 1&0\\0&0\endmatrix\right)
\tag{22} 
$$ 
in $\Bbb B(\frak H^2)$. Then for some  $s$ in $I$ 
we have by the pinching inequality in (iv) that 
$$
Pf(PU^*XUP+s(\bold 1-P))P \le Pf(U^*XU)P= PU^*f(X)UP.
\tag{23}
$$
Since $$U^*XU=\pmatrix \lambda x+(1-\lambda) y 
& (\lambda-\lambda^2)^{1/2}(y-x)\\
(\lambda-\lambda^2)^{1/2}(y-x) 
& \lambda y + (1-\lambda) x\endpmatrix,
\tag{24}
$$
it follows that
$$
\aligned &\pmatrix f(\lambda x + (1-\lambda)y) & 0
\\0 & 0\endpmatrix =Pf(PU^*XUP+s(\bold 1-P))P\\ 
\le\;& PU^*\left(\matrix f(x)&0\\0 & f(y)\endmatrix\right)UP 
=\pmatrix \lambda f(x)+(1-\lambda)f(y) & 0\\0 & 0\endpmatrix. 
\endaligned 
\tag{25} 
$$
\hfill $\square$
\enddemo


\demo{Proof of Corollary 2.3} If $\sum_{k=1}^n 
a^*_ka_k = b \le \bold 1$, put $a_{n+1}=
(\bold 1 -b)^{1/2}$. Then we have a unital 
$(n+1)-$tuple, so with $x_{n+1}=0$ we get
$$
\aligned
&f\left(\sum_{k=1}^n a^*_kx_ka_k\right) 
= f\left(\sum_{k=1}^{n+1} a^*_kx_ka_k\right)
\le\sum_{k=1}^{n+1} a^*_kf(x_k)a_k\\
&\;  =\sum_{k=1}^n a^*_kf(x_k)a_k + a^*_{n+1}f(0)a_{n+1}
\le \sum_{k=1}^n a^*_kf(x_k)a_k.
\endaligned
\tag{26}
$$

Conversely, if (5) is satisfied for all contractive 
$n-$tuples, then -- a fortiori -- it holds for 
unital $n-$tuples, so $f$ is operator convex; and 
with $a=x=0$ we see that $f(0)\le 0\cdot f(0)\cdot 0 =0$.
\hfill$\square$\enddemo


\demo{Proof of Theorem 2.4} Let $ x_k=
\sum_{\sp (x_k)}\lambda E_k(\lambda)$
denote the spectral resolution of $x_k$ for 
$1\le k\le n$. Thus, $E_k(\lambda)$ is the spectral 
projection of $x_k$ on the eigenspace corresponding 
to $\lambda$ if $\lambda$ is an eigenvalue for $x_k$\,; 
otherwise $E_k(\lambda)=0$. For each unit vector $\xi$ 
in $\Bbb C^m$ define the (atomic) probability measure 
$$
\mu_\xi(S)
=\left(\sum_{k=1}^n a_k^* E_k(S)a_k\xi\,\bigg\vert\,\xi\right)
=\sum_{k=1}^n \left(E_k(S)a_k\xi\mid a_k\xi\right)
\tag{27}
$$
for any (Borel) set $S$ in $\Bbb R$. Note now that if  
$y=\sum_{k=1}^n a_k^*x_ka_k$ then
$$
\aligned
&(y\xi|\xi)=\left(\sum_{k=1}^n a_k^* x_k a_k\xi\,\bigg\vert\,\xi\right)\\
=\;&\left(\sum_{k=1}^n\sum_{\sp (x_k)} \lambda E_k(\lambda)
a_k\xi\,\bigg\vert\, a_k\xi\right)=\int\lambda\, d\mu_\xi(\lambda).
\endaligned
\tag{28}
$$
If a unit vector $\xi$ is an eigenvector for $y$, then the
corresponding eigenvalue is $(y\xi|\xi)$, and $\xi$ is also 
an eigenvector for $f(y)$ with correponding eigenvalue 
$(f(y)\xi|\xi)=f((y\xi|\xi))$. In this case we therefore have 
$$
\aligned
&\left(f\left(\sum_{k=1}^n a_k^* x_k
a_k\right)\xi\,\bigg\vert\,\xi\right)=(f(y)\xi|\xi)=f((y\xi|\xi))\\
=\;&f\left(\int\lambda\, d\mu_\xi(\lambda)\right)
\le\int f(\lambda)\, d\mu_\xi(\lambda)\\
=\;&\sum_{k=1}^n\left(\sum_{\sp (x_k)} f(\lambda)E_k(\lambda) a_k\xi
\,\bigg\vert\, a_k\xi\right) =\sum_{k=1}^n\left(a_k^* f(x_k)
a_k\xi\mid\xi\right),
\endaligned
\tag{29}
$$
where we used (28) and the convexity of $f$ -- in form 
of the usual Jensen inequality -- to get the inequality 
in (29). The result in (7) now follows by
summing over an orthonormal basis of eigenvectors for $y$.

Conversely, if (7) holds for some pair of natural numbers 
$n, m$, where $n>1$, then taking $a_i=0$ for 
$i\ge 2$ we see that the inequality holds for $n=2$. Given $s, t$ in
$I$ and $\lambda$ in $[0, 1]$ we define $x=s \bold 1_m$ and
$y=t \bold 1_m$ in $\Bbb M_m(\Bbb C)$. Then with $a=\lambda^{1/2}\bold
1_m$ and $b=(1-\lambda)^{1/2} \bold 1_m$ we get by (7) that 
$$
\aligned
&m f(\lambda x+(1-\lambda)t) 
=\Tr (f(\lambda x+(1-\lambda)t) \bold 1_m)\\ 
=\;& \Tr (f(a^*xa+b^*yb)) \le \Tr (a^*f(x)a+b^*f(y)b) \\
=\; & \Tr ((\lambda f(s)+(1-\lambda)f(t))\bold 1_m)
=m(\lambda f(s)+(1-\lambda)f(t)),
\endaligned
\tag{30}
$$
which shows that $f$ is convex. \hfill$\square$
\enddemo


\subhead{Continuous Fields of Operators}\endsubhead 
Let $\Cal A$ be a $C^*-$algebra of operators on some 
Hilbert space $\frak H$ and $T$ a locally compact Hausdorff 
space. We say that a family $(a_t)_{ t\in T}$ of operators 
in the multiplier algebra $M(\Cal A)$ of $\Cal A$,  
i\.e\. the $C^*-$algebra 
$\{a\in \Bbb B(\frak H)\mid \forall x\in \Cal A \, : \, 
xa+ax\in \Cal A \}$, is a {\it continuous field}, if the function 
$t @>>> a_t$ is norm continuous. If $\mu$ is a Radon measure 
on $T$ and the function $t@>>> \Vert a_t\Vert$ is integrable, 
we can then form the Bochner integral $\int_T a_t\,d\mu(t)$, 
which is the unique element in $M(\Cal A)$ such that  
$$
\varphi \left(\int_T a_t\, d\mu (t)\right) = 
\int_T \varphi (a_t)\, d\mu (t)  \qquad \varphi \in \Cal A^*.
\tag{31}
$$
If all the $a_t$'s belong to $\Cal A$  then also  
$\int_T a_t\,d\mu(t)$ belongs to $\Cal A$. If 
$(a^*_ta_t)_{t\in T}$ is integrable with
integral $\bold 1$ we say that 
$(a_t)_{ t\in T}$ is a {\it unital column field}. 

The transition from sums to continuous fields is 
prompted by the nature of the proof of Theorem 4.1, 
but the interested reader can easily verify that 
also Theorem 2.1 is valid for continuous fields. 
We finally note that the restriction to continuous 
fields is handy, but not necessary. In \cite{{\bf 10}} 
we shall generalize the setting to arbitrary weak* 
measurable fields. 


\subhead{Centralizers}\endsubhead Recall that the 
{\it centralizer} of a positive functional 
$\varphi$ on a $C^*-$al\-ge\-bra $\Cal A$ is the closed
$^*-$subspace $\Cal A^{\varphi} = 
\{y\in \Cal A\mid \forall x\in \Cal A\,
:\, \varphi(xy) =\varphi(yx)\}$.  In general this is 
not an algebra, but if $y_1,\dots, y_n$ are pairwise 
commuting, self-adjoint elements in $\Cal A^{\varphi}$ 
then the $C^*-$algebra they generate is contained in 
$\Cal A^{\varphi}$. Evidently the size of 
$\Cal A^{\varphi}$ measures the extent to which 
$\varphi$ is a trace. The fact we shall utilize is 
that even if an element $x$ is outside  
$\Cal A^{\varphi}$ the functional will behave 
"trace-like" on the subspace spanned by
$\Cal A^{\varphi}x\Cal A^{\varphi}$.  

If $\varphi$ is unbounded, but lower semi-continuous 
on $\Cal A_+$ and finite on the minimal 
dense ideal $K(\Cal A)$ of $\Cal A$, we 
define $\Cal A^{\varphi}=\{y\in\Cal A
\mid\forall x\in K(\Cal A)\,:\,\varphi(xy)=\varphi(yx)\}$.

\bigskip

\proclaim{Theorem 4.1} Let $(x_t)_{t\in T}$ be a 
bounded, continuous field on a locally compact 
Hausdorff space $T$ consisting of self-adjoint elements
in a $C^*-$algebra $\Cal A$ with $\sp (x_t) \subset I$.  
Furthermore, let $(a_t)_{t\in T}$ be a unital column 
field in $M(\Cal A)$ with respect to some Radon measure 
$\mu$ on $T$. Then for each continuous, convex function
$f$ defined on $I$ and every positive functional 
$\varphi$ that contains the element 
$y=\int_T a^*_tx_ta_t\,d\mu(t)$ in its centralizer 
$\Cal A^{\varphi}$, i\.e\. $\varphi (xy) = \varphi (yx)$ 
for all $x$ in $\Cal A$, we have the inequality:
$$
\varphi\left(f\left(\int_T a_t^*x_t a_t\,d\mu (t)\right)\right)
\le \varphi\left(\int_T a^*_t f(x_t)a_t\,d\mu (t)\right).
\tag{32}
$$

If $\varphi$ is unbounded, but lower semi-continuous on 
$\Cal A_+$ and finite on the minimal dense ideal $K(\Cal A)$ 
of $\Cal A$, the result still holds if $f\ge 0$, even 
though the function may now attain infinite values. 
\endproclaim

\demo{Proof} Let $\Cal C = C_o (S)$ denote the commutative 
$C^*-$subalgebra of $\Cal A$ generated by $y$, 
and let $\mu_\varphi$ be the finite Radon measure on 
the locally compact Hausdorff space $S$ defined, via the 
Riesz representation theorem, by 
$$
\int_S z(s)\,d\mu_{\varphi} (s)=\varphi(z)\qquad z\in \Cal C=C_o(S).
\tag{33}
$$
Since for all $(x, z)$ in $M(\Cal A)_+\times \Cal C_+$ we have 
$\varphi(xz)=\varphi(z^{1/2}xz^{1/2})$ it follows that
$$
0\le \varphi(xz) \le \Vert x \Vert \varphi(z).
\tag{34}
$$
Consequently the functional $z\to\varphi(xz)$ on $\Cal C$ 
defines a Radon measure on $S$ dominated by a multible 
of $\mu_{\varphi}$, hence determined by a unique element 
$\Phi(x)$ in $L^\infty_{\mu_\varphi} (S)$. By linearization
this defines a conditional expectation 
$\varphi\colon M(\Cal A)\to L^\infty_{\mu_\varphi}(S)$ 
(i\.e\. a positive, unital module map) such that 
$$
\int_S z(s)\Phi(x)(s)\,d\mu_\varphi (s)
=\varphi(zx),\qquad z\in \Cal C
\quad x\in M(\Cal A).
\tag{35}
$$
Inherent in this formulation is the fact that if 
$z\in \Cal C = C_o(S)$, then $\Phi(z)$ is the 
natural image of $z$ in $L^\infty_{\mu_\varphi}(S)$.
In particular, $z(s) = \Phi(z)(s)$ for almost all $s$ in $S$. 

Observe now that since the $C^*-$algebra $C_o(I)$ is 
separable, we can for almost every $s$ in $S$ 
define a Radon measure $\mu_s$ on $I$ by
$$
\int_I g(\lambda)\, d\mu_s (\lambda) 
= \Phi\left( \int_T a^*_t g(x_t)a_t\,d\mu(t)\right)(s) \qquad g\in C(I).
\tag{36}
$$
As $\int_Ta^*_ta_t\,d\mu(t)=\bold 1$, this is actually 
a probability measure. If we take $g(\lambda)= \lambda$, then 
$$
\int_I \lambda\,d\mu_s(\lambda) = \Phi\left(\int_T a^*_t x_t a_t\,
d\mu(t)\right)(s) = \Phi(y) (s) = y(s).
\tag{37}
$$
Since $y\in C$ we get by (37) and (36) -- using also
the convexity of $f$ in form of the standard Jensen 
inequality -- that
$$
\aligned
&f(y)(s)=f(y(s))=f\left(\int_I \lambda\,d\mu_s(\lambda)\right) \\ 
\le\;& \int_I f(\lambda)\,d\mu_s(\lambda) = 
\Phi\left(\int_T a^*_tf(x_t)a_t\,d\mu(t)\right)(s).
\endaligned
\tag{38}
$$
Integrating over $s$, using (35), now gives the desired result:
$$
\aligned
\varphi(f(y)) =\;& \int_S f(y)(s)\, d\mu_\varphi(s)\\
\le\;&\int_S \Phi\left(\int_T  a^*_t f(x_t)a_t\, 
d\mu(t)\right) (s)\, d\mu_\varphi(s) \\
=\;&\int_T\int_S \Phi \left( a^*_tf(x_t)a_t \right)(s)\, d\mu_\varphi
(s)\,d\mu(t)\\ =\;&\int_T \varphi\left( a^*_tf(x_t)a_t\right)\, d\mu(t)=
\varphi\left( \int_T a^*_t f(x_t)a_t\, d\mu(t)\right). 
\endaligned 
\tag{39}
$$
 
Having proved the finite case, let us now assume 
that $\varphi$ is unbounded, but lower semi-continuous 
on $\Cal A_+$ and finite on the minimal dense ideal 
$K(\Cal A)$ of $\Cal A$. This -- by definition -- means that 
$\varphi(x)<\infty$ if $x\in \Cal A_+$ and $x=xe$ for 
some $e$ in $\Cal A_+$, because $K(\Cal A)$ is the 
hereditary $^*-$subalgebra of $\Cal A$ generated by such
elements, cf\. \cite{{\bf 17}, 5.6.1}. 
Restricting $\varphi$ to $\Cal C$ we therefore 
obtain a unique Radon measure $\mu_\varphi$ on 
$S$ such that    
$$
\int_S z(t)\,d\mu_\varphi (t) = \varphi (z) \qquad y\in \Cal C.
\tag{40}
$$
Inspection of the proof above now shows that the 
Jensen trace inequality still holds if only $f \ge 0$, 
even though $\infty$ may now occur in the inequality. 
\hfill$\square$
\enddemo

\bigskip

\demo{Proof of Theorem 2.7} Evidently this (like Theorem 2.4) 
is a special case of Theorem 4.1, where the continuous field is 
replaced by a finite sum and the functional $\varphi$ is a trace, 
so that $\Cal A^{\varphi}=\Cal A$. \hfill$\square$
\enddemo

\newpage

\enddocument